\documentclass[12pt,a4paper]{amsart}
\usepackage[utf8]{inputenc}	
\usepackage[T1]{fontenc}
\usepackage[in,headings]{fullpage}
\usepackage{amsbsy, amscd, amsfonts, amsmath, amsrefs, amssymb, amsthm}
\usepackage{graphicx}
\usepackage{float}
\usepackage{color}   
\usepackage[colorlinks=true,linkcolor=blue,citecolor=blue,urlcolor=blue]{hyperref}
\usepackage{comment}
\excludecomment{A}

\newtheorem{theorem}{Theorem}

\newtheorem{corollary}[theorem]{Corollary}
\newtheorem{proposition}[theorem]{Proposition}
\newtheorem{lemma}[theorem]{Lemma}

\theoremstyle{definition}
\newtheorem{definition}[theorem]{Definition}
\newtheorem{remark}[theorem]{Remark}

\theoremstyle{remark}

\newcommand{\C}{\mathbf{C}}
\newcommand{\Z}{\mathbf{Z}}

\newcommand{\N}{\mathbf{N}}

\renewcommand{\Re}{\mathop{\mathrm{Re}}\nolimits}
\renewcommand{\Im}{\mathop{\mathrm{Im}}\nolimits}

\newcommand{\Lzeta}{\mathop{\mathcal L }\nolimits}
\newfont{\cmbsy}{cmbsy10}
\newfont{\cmmib}{cmmib10}

\begin{document}

\title[An expression for Riemann Siegel function]
{An expression for Riemann Siegel function}
\author[Arias de Reyna]{J. Arias de Reyna}
\address{%
Universidad de Sevilla \\ 
Facultad de Matem\'aticas \\ 
c/Tarfia, sn \\ 
41012-Sevilla \\ 
Spain.} 




\def\checked{}

\email{arias@us.es, ariasdereyna1947@gmail.com}


\begin{abstract}
There are many analytic functions $U(t)$ satisfying $Z(t)=2\Re\bigl\{ e^{i\vartheta(t)}U(t)\bigr\}$. Here, we consider an entire function $\Lzeta(s)$
such that $U(t)=\Lzeta(\frac12+it)$ is one of the simplest among them. 
We obtain an expression for the Riemann-Siegel function $Z(t)$ in terms of the zeros of $\Lzeta(s)$. Implicitly, the function $\Lzeta(s)$ is considered by Riemann in \cite{R}.

Riemann spoke of having used an expression for $\Xi(t)$ in his demonstration that most of the non-trivial zeros of the zeta function lie on the critical line. Therefore, any expression deserves a study. 
\end{abstract}

\maketitle

\section{Introduction.}

Riemann stated \cite{R2}*{Letter to Weierstrass} that he has proved that the number of zeta zeros in the critical line is about $\frac{T}{2\pi}\log\frac{T}{2\pi}-\frac{T}{2\pi}$. He also said that the proof follows from a new representation of the function $\Xi(t)$ that he had not simplified  sufficiently to be able to publish it. In this paper, we develop a representation that follows easily from results contained in Riemann's paper \cite{R}. 

These statements by Riemann have always been received with much skepticism. See, for example,  Edwards \cite{Edwards}*{p.~164}. But after reading Siegel \cite{Siegel} and writing \cite{A49}, I do not believe any more in a naive Riemann.

This representation is related to the expansion due to Lavrik. Therefore, some of our results are well known. We see that the end representation of the Riemann-Siegel function $Z(t)$ does not appear very useful for computation. But because of the possible connections with Riemann's cryptic statements they deserve to be studied.

\section{New expression for \texorpdfstring{$\Xi(t)$}{X(t)}.}

In Riemann's paper  can be found the expression
\[\checked{\pi^{-s/2}\Gamma(s/2)\zeta(s)=\frac{1}{s(s-1)}+\int_1^{+\infty}\psi(x)\bigl(x^{\frac{s}{2}-1}+x^{-\frac{1+s}{2}}\bigr)\,dx.}\]
Following the reasoning of Riemann, Lavrik  extended this.
Put 
\[\theta(z)=\sum_{n\in\Z}e^{-\pi n^2z},\qquad\Re z>0.\]
This function is essentially the ordinary theta function. It is analytic in the  half-plane $\Re z>0$. And we have $\psi(x)=(\theta(x)-1)/2$. 

\begin{theorem}
For any complex number $\tau$ with $\Re(\tau)>0$ and any $s\notin\{0,1\}$ we have
\begin{equation}\label{E:lavrik}
\checked{\begin{aligned}
\pi^{-s/2}\Gamma(\tfrac{s}{2})\zeta(s)&=
-\frac{\tau^{s/2}}{s}+\tau^{s/2}\int_1^\infty x^{s/2}\frac{\theta(\tau x)-1}{2}\frac{dx}{x}\\
&-\frac{\tau^{-(1-s)/2}}{1-s}+
\tau^{-(1-s)/2}\int_{1}^\infty 
y^{(1-s)/2}\frac{\theta(y/\tau)-1}{2}\frac{dy}{y}.
\end{aligned}}
\end{equation}
\end{theorem}

\begin{proof}
In Riemann's paper we find the expression
\[\pi^{-s/2}\Gamma(\tfrac{s}{2})\zeta(s)=\int_0^\infty x^{s/2}\frac{\theta(x)-1}{2}\frac{dx}{x},\qquad \sigma>1.\]
Given any number $\tau$ with $\Re(\tau)>0$ we may rotate  the line of integration 
to make it pass through  $\tau$.
We obtain
\begin{equation}\label{E:partial}
\int_0^\infty x^{s/2}\frac{\theta(x)-1}{2}\frac{dx}{x}=\int_0^\tau x^{s/2}\frac{\theta(x)-1}{2}\frac{dx}{x}+\int_\tau^\infty x^{s/2}\frac{\theta(x)-1}{2}\frac{dx}{x}.
\end{equation}

In the first integral, we may integrate a part and then 
change the variables putting $x =\tau/y$ 
\[\int_0^\tau x^{s/2}\frac{\theta(x)-1}{2}\frac{dx}{x}=-\frac{\tau^{s/2}}{s}+\frac12\int_{1}^\infty 
(\tfrac{\tau}{y})^{s/2}\theta(\tau/y)\frac{dy}{y}\]
By the functional equation $\theta(\tau/y)=\sqrt{y/\tau}\theta(y/\tau)$. Therefore,
\[\int_0^\tau x^{s/2}\frac{\theta(x)-1}{2}\frac{dx}{x}=-\frac{\tau^{s/2}}{s}+\frac12\int_{1}^\infty 
(\tfrac{y}{\tau})^{(1-s)/2}\theta(y/\tau)\frac{dy}{y}.\]
Since $\Re(1-\sigma)<0$,  we may insert  the $1/2$ term again
\begin{multline*}
-\frac{\tau^{s/2}}{s}+\frac12\int_{1}^\infty 
(\tfrac{y}{\tau})^{(1-s)/2}\theta(y/\tau)\frac{dy}{y}=
-\frac{\tau^{s/2}}{s}-\frac{\tau^{-(1-s)/2}}{1-s}+\int_{1}^\infty 
(\tfrac{y}{\tau})^{(1-s)/2}\frac{\theta(y/\tau)-1}{2}\frac{dy}{y}
\end{multline*}
Substituting into \eqref{E:partial} and changing variables in one of the integrals, we obtain \eqref{E:lavrik}.
Since the two integrals are entire functions of $s$, the equality extends to 
any $s$ not equal to $0$ or $1$. 
\end{proof}

\begin{corollary}\label{C:cor}
For any complex number $\tau$ with $\Re(\tau)>0$ and any $s\notin\{0,1\}$ we have
\begin{equation}\label{E:firstexpr}
\checked{\begin{aligned}
\pi^{-s/2}\Gamma(\tfrac{s}{2})\zeta(s)&=
-\frac{\tau^{s/2}}{s}+\pi^{-s/2}\sum_{n=1}^\infty 
\frac{\Gamma(\frac{s}{2},\pi n^2\tau)}{n^s}
\\
&-\frac{\tau^{-(1-s)/2}}{1-s}+
\pi^{-(1-s)/2}\sum_{n=1}^\infty 
\frac{\Gamma(\frac{1-s}{2},\pi n^2/\tau)}{n^{1-s}}.
\end{aligned}}
\end{equation}
\end{corollary}
\begin{proof}
For $\Re\tau>0$, it is easy to justify the integration term by term to obtain
\[\int_1^\infty x^{s/2}\frac{\theta(\tau x)-1}{2}\frac{dx}{x}=
\sum_{n=1}^\infty\int_1^\infty x^{s/2}e^{-\pi n^2 \tau x}\frac{dx}{x}=
\frac{\Gamma(s/2,\pi n^2\tau)}{(\pi\tau)^{s/2}n^s}.\]
We can transform in an analogous way the second part of \eqref{E:lavrik} to get \eqref{E:firstexpr}.
\end{proof}

\begin{definition}
For any $\tau$ with $\Re(\tau)>0$ and $s\in\C$ we define $\Lzeta(\tau,s)$
by the equation 
\begin{equation}\label{E:Uno}
\pi^{-s/2}\Gamma(\tfrac{s}{2})\Lzeta(\tau,s)=
-\frac{\tau^{s/2}}{s}+\tau^{s/2}\int_1^\infty x^{s/2}\frac{\theta(\tau x)-1}{2}\frac{dx}{x}.
\end{equation}
We are interested specially in the case $\tau=1$ in which case we will put 
$\Lzeta(s)=\Lzeta(1,s)$.
\end{definition}
It is easily seen that the integral is an entire function of $s$.
Thus, \[\Lambda(\tau,s):=\pi^{-s/2}\Gamma(s/2)\Lzeta(\tau,s)\] is a meromorphic
function with a unique pole at the point $s=0$. This pole is simple
and with residue $-1$. The above equation \eqref{E:Uno} also shows
that $\Lzeta(\tau,-2n)=0$ for $n\in \N$. This we call the
trivial zeros of $\Lzeta(\tau,s)$. It also follows that $\Lzeta(\tau,s)$ is an
entire function.

With this notation, Lavrik equation \eqref{E:lavrik} can be written as 
\begin{equation}\label{E:InRiemann}
\checked{\pi^{-s/2}\Gamma(\tfrac{s}{2})\zeta(s)=\pi^{-s/2}\Gamma(\tfrac{s}{2})\Lzeta(\tau,s)+\pi^{-(1-s)/2}\Gamma(\tfrac{1-s}{2})\Lzeta(1/\tau,1-s).}
\end{equation}
And Corollary \ref{C:cor} or his proof implies that 
\begin{equation}\label{E:seriesrep}
\Lzeta(\tau,s)=-\frac{(\pi\tau)^{s/2}}{s\;\Gamma(\frac{s}{2})}+\sum_{n=1}^\infty \frac{\Gamma(\frac{s}{2},\pi n^2\tau)}{\Gamma(\frac{s}{2})}\frac{1}{n^s}.
\end{equation}

As an example of Mellin-Barnes integrals \href{http://dlmf.nist.gov/8.6.iii}{see DLMF}  we have 
\begin{lemma}
For $\Re z>0$ and $a\ne 0$, $-1$, $-2$, \dots, we have
\begin{equation}\label{E:mellinbarnes}
\checked{\Gamma(a,z)=
\frac{1}{2\pi i}\int_{c-i\infty}^{c+i\infty}\Gamma(s+a)\frac{z^{-s}}{s}\,ds}
\end{equation}
where $c$ is a real constant and the path of integration is indented (if necessary) so that it  is to the right of all poles. 
\end{lemma}

\begin{proposition}
For $s\ne 0$ we have
\begin{equation}\label{E:neweq}
\checked{\pi^{-\frac{s}{2}}\Gamma({\textstyle
\frac{s}{2}})\Lzeta(\tau,s)=-\frac{\tau^{s/2}}{s}+\frac{\tau^{s/2}}{2\pi i}\int_{c-i\infty}^{c+i\infty}\pi^{-\frac{z}{2}}\Gamma(\tfrac{z}{2})\zeta(z)\frac{\tau^{-z/2}}{z-s}\,dz,}
\end{equation}
where the integration path is a vertical line $\Re z= c$ with $c>\max(\sigma,1)$.
\end{proposition}
\begin{proof}
By \eqref{E:seriesrep} we have 
\[\pi^{-s/2}\Gamma(\tfrac{s}{2})\Lzeta(\tau,s)=-\frac{\tau^{s/2}}{s}+\pi^{-s/2}\sum_{n=1}^\infty \frac{\Gamma(\frac{s}{2},\pi n^2\tau)}{n^s}.\]
By the Mellin Barnes equation \eqref{E:mellinbarnes} 
\[\Gamma(\tfrac{s}{2},\pi n^2\tau)=\frac{1}{2\pi i}\int_{c-i\infty}^{c+i\infty}\Gamma(z+\tfrac{s}{2})\frac{(\pi n^2\tau)^{-z}}{z}\,dz,\]
where $c$ is a real constant such that $c\ge \max(-\sigma/2,0)$. Changing the variables $z\mapsto \frac{z-s}{2}$ 
\[\Gamma(\tfrac{s}{2},\pi n^2\tau)=\frac{1}{2\pi i}\int_{c-i\infty}^{c+i\infty}\Gamma(\tfrac{z}{2})\frac{(\pi n^2\tau)^{-\frac{z-s}{2}}}{z-s}\,dz,\]
where now $c>\max(\sigma,0)$.  Therefore, 
\begin{align*}
\pi^{-s/2}\sum_{n=1}^\infty \frac{\Gamma(\frac{s}{2},\pi n^2\tau)}{n^s}&=
\frac{\pi^{-s/2}}{2\pi i}\sum_{n=1}^\infty \frac{1}{n^s}\int_{c-i\infty}^{c+i\infty}\Gamma(\tfrac{z}{2})\frac{(\pi n^2\tau)^{-\frac{z-s}{2}}}{z-s}\,dz\\
&=\frac{\tau^{s/2}}{2\pi i}\sum_{n=1}^\infty \int_{c-i\infty}^{c+i\infty}\pi^{-\frac{z}{2}}\Gamma(\tfrac{z}{2})\frac{1}{n^z}\frac{\tau^{-\frac{z}{2}}}{z-s}\,dz\\
\end{align*}
Taking $c>1$ we may interchange sum and integral to get 
\[\pi^{-s/2}\sum_{n=1}^\infty \frac{\Gamma(\frac{s}{2},\pi n^2\tau)}{n^s}=\frac{\tau^{s/2}}{2\pi i}\int_{c-i\infty}^{c+i\infty}\pi^{-\frac{z}{2}}\Gamma(\tfrac{z}{2})\zeta(z)\frac{\tau^{-z/2}}{z-s}\,dz.\]
This proves \eqref{E:neweq}.
\end{proof}

With the usual notation, we have
\begin{align*}
\xi(s)&=\frac{s(s-1)}{2}\pi^{-s/2}\Gamma(s/2)\zeta(s)
\\&=\frac{s(s-1)}{2}\pi^{-s/2}\Gamma(\tfrac{s}{2})\Lzeta(s)+
\frac{s(s-1)}{2}\pi^{-(1-s)/2}\Gamma(\tfrac{s}{2})\Lzeta(1-s),
\end{align*}
and for $z\in\C$ we define $\Xi(z)=\xi(\frac12+iz)$

Defining for any complex $z$ 
\[F(z):=\frac{s(s-1)}{2}\pi^{-s/2}\Gamma(\tfrac{s}{2})\Lzeta(s),\qquad s=\tfrac12+iz,\]
we get 
\[\Xi(z)=F(z)+F(-z).\]
Here, $F(z)$ is an entire function.

\begin{corollary}
For real $t$, we have
\begin{equation}\label{E:ZL}
\checked{Z(t)=2\Re\bigl\{ e^{i\vartheta(t)}\Lzeta(\tfrac12+it)\bigr\}.}
\end{equation}
\end{corollary}

\begin{proof}
For $s=\tfrac12+it$ we have
$\pi^{-s/2}\Gamma(s/2)=f(t)e^{i\vartheta(t)}$ with $f(t)>0$. Thus,
\eqref{E:InRiemann} can be written as
\begin{displaymath}
f(t)e^{i\vartheta(t)}\zeta(\tfrac12+it)=
f(t)e^{i\vartheta(t)}\Lzeta(\tfrac12+it)+f(t)e^{-i\vartheta(t)}\Lzeta(\tfrac12-it).
\end{displaymath}
We simplify $f(t)$,  then the two terms are complex conjugate and
its sum is twice the real part.
\end{proof}

\begin{theorem}\label{T:product}
$s\pi^{-s/2}\Gamma(s/2)\Lzeta(s)$ is an entire function of order $1$
that has a unique real zero $b_0$ and an infinite set of complex
conjugate zeros $b_n$ and $b_{-n}$. We have the following expression
as an infinite product
\begin{equation}
\Lambda(s)=\pi^{-s/2}\Gamma(s/2)\Lzeta(s)=-\frac{e^{-\alpha
s}}{s}\prod_{n\in\Z} \Bigl(1-\frac{s}{b_n}\Bigr) e^{s/b_n}.
\end{equation}
\end{theorem}

\begin{proof}
Put $\Lambda(s)=\pi^{-s/2}\Gamma(s/2)\Lzeta(s)$.  By \eqref{E:Uno} with $\tau=1$,
we get a bound of the function $s\Lambda(s)$ in the following way
(with $|s|=r$)
\begin{align*}
|s\pi^{-s/2}\Gamma(s/2)\Lzeta(s)|&\le1+\Bigl|s\int_1^{+\infty}y^{s/2}\Bigl(\sum_{n=1}^\infty
e^{-\pi n^2 y}\Bigr)\frac{dy}{y}\Bigr|\\
&\le1+2r\int_1^{+\infty} y^{r/2}e^{-\pi y}\frac{dy}{y}\\
&\le 1+2r\pi^{-r/2}\int_0^{+\infty}x^{r/2}e^{-x} \frac{dx}{x}\le
1+2r\pi^{-r/2} \Gamma(r/2).
\end{align*}
 It is clear that this bound is optimal since for real and positive
$s$  it is almost an equality. Therefore, if $M(r)$ denotes the maximum of
$|s\pi^{-s/2}\Gamma(s/2)\Lzeta(s)|$ for $|s|\le r$ we will have
$M(r)>cr\pi^{-r/2}\Gamma(r/2)$.  From all this it is clear that
$s\Lambda(s)=s\pi^{-s/2}\Gamma(s/2)\Lzeta(s)$ has order $1$.

The function $s\Lambda(s)$  has an infinite number of zeros, because
if not we will have $s\Lambda(s)=P(s) e^{Q(s)}$ for some polynomials
$P$ and $Q$, $Q$ will be of degree one, and this contradicts that
$M(r)>cr\pi^{-r/2}\Gamma(r/2)$.

From \eqref{E:Uno} we easily find that $\Lambda(x)$ is an
increasing function for real $x$.
\begin{displaymath}
\lim_{x\to0^-}\Lambda(x)=+\infty;\qquad
\lim_{x\to0^+}\Lambda(x)=-\infty.
\end{displaymath}
Since it is positive for $x<0$, it follows that at most it has a
real zero. This zero exists because
$\lim_{x\to+\infty}\Lambda(x)=+\infty$. Numerically, it is easy to
compute that this zero is approximately $b_0\approx11.25170908146$.
The other zeros can be paired as complex conjugate. Hence, we can
number the zeros as $b_n$ with $b_{-n}=\overline{b_n}$. We assume
that they are numbered in such a way that $|b_n|\le|b_{n+1}|$ for
$n\ge1$.

Now it is clear that we can put
\begin{displaymath}
s\Lambda(s)=Ae^{-\alpha s}\prod_{n\in\Z} \Bigl(1-\frac{s}{b_n}\Bigr)
e^{s/b_n}.
\end{displaymath}
Since
\begin{displaymath}
\lim_{s\to0} s\Lambda(s)=-1+\lim_{s\to0} s\int_1^{+\infty}
y^{s/2}\Bigl(\sum_{n=1}^\infty e^{-\pi n^2 y}\Bigr)\frac{dy}{y}=-1.
\end{displaymath}
we know that $A=-1$.

By taking logarithms
\begin{displaymath}
\log\{-s\Lambda(s)\}=-\alpha s-\sum_{k=2}^\infty\Bigl(\sum_{n\in\Z}
\frac{1}{b_n^k}\Bigr)\frac{s^k}{k}.
\end{displaymath}
On the other hand, by  \eqref{E:Uno} for  small enough $s$ 
\begin{displaymath}
\log\{-s\Lambda(s)\}=-\sum_{k=1}^\infty\Bigl(\int_1^{+\infty}
y^{s/2}\Bigl(\sum_{n=1}^\infty e^{-\pi n^2
y}\Bigr)\frac{dy}{y}\Bigr)^k\frac{s^k}{k}
\end{displaymath}
\begin{displaymath}
=-\sum_{k=1}^\infty\Bigl\{\sum_{m=0}^\infty \Bigl(\int_1^{+\infty}
(\log y)^m\Bigl(\sum_{n=1}^\infty e^{-\pi n^2
y}\Bigr)\frac{dy}{y}\Bigr)\frac{1}{m!}\Bigl(\frac{s}{2}\Bigr)^m\Bigr\}^k
\frac{s^k}{k}.
\end{displaymath}
Thus  $-\alpha$ is the coefficient of  $s$ in this expansion. An
easy computation gives us
\begin{displaymath}
\alpha=\int_1^{+\infty}\Bigl(\sum_{n=1}^\infty e^{-\pi n^2
y}\Bigr)\frac{dy}{y}\approx0.010906559198968892180277118987156\qedhere
\end{displaymath}
\end{proof}

\begin{figure}[H]
\includegraphics*[width=7cm]{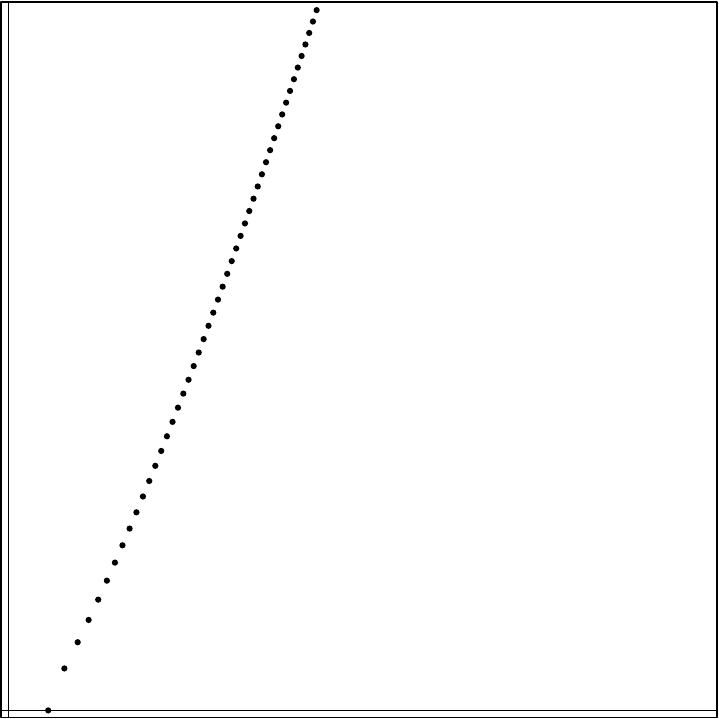}\\
\caption{Zeros of $s\Lambda(s)$ on $(-2,200)\times(-2,200)$}
\label{F:Uno}
\end{figure}

It is not difficult to compute the zeros of $\Lzeta$ numerically.
The zeros of $\Lambda(s)$ appear to be contained in the first and fourth
quadrants. We can see in  Figure \ref{F:Uno} those in the first
quadrant. The others are the complex conjugates of these.

We can see that apparently the zeros $(b_n)$  are neatly arranged.
We have calculated 807 zeros,   those with $\Im b_n\le 1972$. The list of these 807 zeros can be found at the end of this \texttt{57-Lavrik-v4.tex} file. Based on these scant data, the
number of zeros $N(x)$ with $|b_n|\le x$ and $n\ge0$ appears to be
\begin{displaymath}
N(x)=\frac{x}{4\pi}\log\frac{x}{2\pi}-\frac{x}{4\pi}+
\sqrt{\frac{x}{4\pi}\log\frac{x}{2\pi}-\frac{x}{4\pi}}+R,
\end{displaymath}
where $R$ is small in the range considered. 
It also appears that if $b=\beta+i \gamma$, then
\begin{displaymath}
\gamma\approx \beta\log\frac{\beta}{4\pi} +R;\qquad
\frac{\gamma}{\log(\gamma/2\pi)}\ge\frac{2\beta}{\pi},
\end{displaymath}
where the error $R$ is relatively small.

\begin{proposition}\label{P:locZeros}
The zeros of the function $s\pi^{-s/2}\Gamma(s/2)\Lzeta(s)$ are
contained in the half-plane $\Re s\ge b_0\approx11.25$.
\end{proposition}

\begin{proof}
By \eqref{E:Uno}
\begin{displaymath}
s\Lambda(s)=s\pi^{-s/2}\Gamma(s/2)\Lzeta(s)=-1+\sum_{n=1}^\infty
s\int_1^{+\infty} y^{s/2}e^{-\pi n^2 y}\frac{dy}{y}.
\end{displaymath}
Integrating by parts, we get
\begin{displaymath}
s\int_1^{+\infty} y^{s/2}e^{-\pi n^2 y}\frac{dy}{y}=-2e^{-\pi
n^2}+2\pi n^2\int_1^{+\infty}y^{s/2}e^{-\pi n^2 y}\,dy.
\end{displaymath}
Thus, if $A=1+2\sum_{n=1}^\infty e^{-\pi n^2}$ we have
\begin{displaymath}
\bigl|s\Lambda(s)+A\bigr|\le\sum_{n=1}^\infty 2\pi
n^2\int_1^{+\infty}y^{\sigma/2}e^{-\pi n^2 y}\,dy.
\end{displaymath}
Hence, if $s$ is a zero of the function, we will have
\begin{displaymath}
A\le \sum_{n=1}^\infty 2\pi n^2\int_1^{+\infty}y^{\sigma/2}e^{-\pi
n^2 y}\,dy.
\end{displaymath}
Now, on the left, we have an increasing function of $\sigma$. This
function increases from $0$ for $\sigma=-\infty$ to $+\infty$ for
$\sigma=+\infty$. Hence  $\sigma\ge\sigma_0$ where $\sigma_0$ is the
value of $\sigma$ for which  we have equality  in the above
inequality. Since
\begin{displaymath}
A= \sum_{n=1}^\infty 2\pi n^2\int_1^{+\infty}y^{\sigma_0/2}e^{-\pi
n^2 y}\,dy
\end{displaymath}
imply that $\Lambda(\sigma_0)=0$. We have that $\sigma_0=b_0$ the
only real zero of $\Lambda(s)$.
\end{proof}

\begin{figure}[H]
  \includegraphics[width=0.49\hsize]{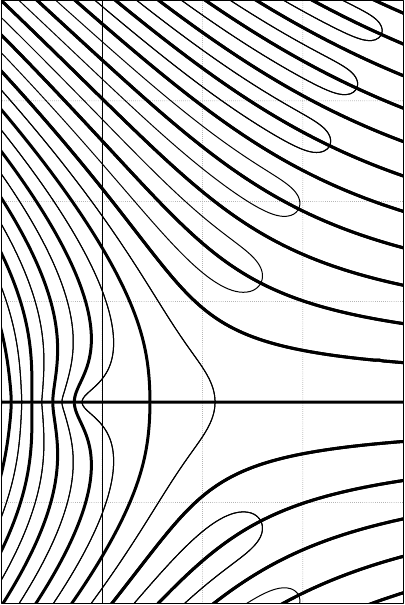}
  \includegraphics[width=0.49\hsize]{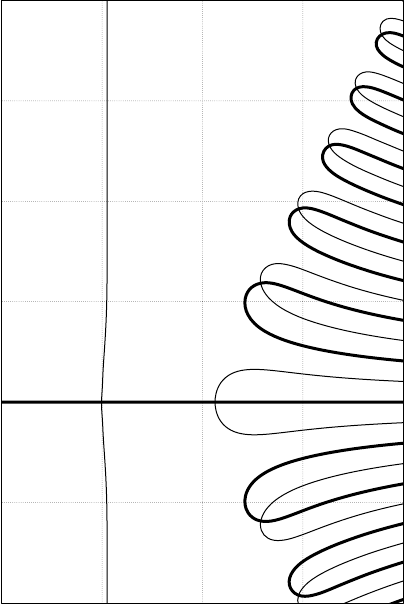}\\
  \caption{$\Lzeta(s)$ and $\pi^{-s/2}\Gamma(s/2)\Lzeta(s)$ on $(-10,30)\times(-20,40)$}
  \label{F:Xray}
\end{figure}

To give an idea of the behavior of the functions 
$\pi^{-s/2}\Gamma(s/2)\Lzeta(s)$ and  $\Lzeta(s)$ we represent in Figure \ref{F:Xray} the X-ray of these functions. On the thick lines  these functions are real,  and on the thin lines they are imaginary.

It is interesting to note that on the graph of the function $\pi^{-s/2}\Gamma(s/2)\Lzeta(s)$ on the left, we see an imaginary line that almost coincides with the critical
line $\sigma=\tfrac12$.  

This imaginary line is sometimes located to the right and sometimes to the left of the critical line $\sigma=\frac12$. At the points where this imaginary line cuts the critical line, we have $Z(t)=0$, that is $\zeta(\frac12+it)=0$.

The function is monotonous on the lines. It follows that on
the critical line the function
$\Lambda(s)=\pi^{-s/2}\Gamma(s/2)\Lzeta(s)$ is almost purely
imaginary and that \[\Im\bigl\{\pi^{-s/2}\Gamma(s/2)\Lzeta(s)\bigr\}\]
is a monotonous decreasing function for $t> t_0$.  Typical values
\begin{align*}
\Lambda(1/2+100i)&\approx  6.844655\;10^{-16} + i\; 0.010862286;\\
\Lambda(1/2+1000i)&\approx 7.955229 \;10^{-17}  + i\; 0.0010864328.
\end{align*}

Since $\Lambda(s)$ does not vanish on the critical line, we can
define $\arg\Lambda(\tfrac12+it)$ as a continuous (real analytic)
function of the real variable $t$. We know that $\Lambda(\tfrac12+it)$
is almost purely imaginary, therefore, we can fix this function
$\arg\Lambda(\tfrac12+it)$ in such a way that its value is almost
$\pi/2$ for big $t>0$.  We prove this in the following

\begin{proposition}\label{P:argf}
There exists a continuous determination of $\arg\Lambda(\tfrac12+it)$
such that for $t>10$ we have
\begin{displaymath}
\Bigl|\arg \Lambda(\tfrac12+it)-\frac{\pi}{2}\Bigr|\le
\arcsin\frac{0.2}{A}=0.19.
\end{displaymath}
\end{proposition}

\begin{proof}
With the notations of the proof of Proposition \ref{P:locZeros} we
have
\begin{displaymath}
s\Lambda(s)=-A+\sum_{n=1}^\infty 2\pi n^2\int_1^{+\infty}
y^{s/2}e^{-\pi n^2y}\,dy.
\end{displaymath}
Thus,
\begin{displaymath}
|s\Lambda(s)+A|\le \sum_{n=1}^\infty 2\pi n^2\int_1^{+\infty}
y^{\sigma/2}e^{-\pi n^2y}\,dy.
\end{displaymath}
For $\sigma=1/2$, we get  that for $s=\tfrac12+it$, we have
\begin{displaymath}
|s\Lambda(s)+A|\le 0.0922;\qquad  \Bigl|\Lambda(s)+\frac{A}{s}\Bigr|\le
\frac{0.0922}{|s|}.
\end{displaymath}
It follows that for $t>10$
\begin{align*}
\Bigl|\Lambda(s)-\frac{iA}{t}\Bigr|&\le
\Bigl|\Lambda(s)+\frac{A}{\tfrac12+it}\Bigr|+\Bigl|\frac{A}{\tfrac12+it}+\frac{iA}{t}\Bigr|\\
&\le\frac{0.0922}{t}+\frac{A}{2t^2}\le \frac{0.2}{t}
\end{align*}

Thus, for $t>10$ there exists a determination of
$\arg\Lambda(\tfrac12+it)$ such that
\begin{displaymath}
\Bigl|\arg \Lambda(\tfrac12+it)-\frac{\pi}{2}\Bigr|\le
\arcsin\frac{0.2}{A}=0.185144288\cdots
\end{displaymath}
since the constant $A=1.0864348112\dots $. It is also clear that
this argument is a continuous function, because it can be seen as the
composition of the usual $\arg$ function, defined on the complex
plane except the negative real axis, with the function $\Lambda$
\end{proof}

We compute  $\Lambda(1/2)=-1.988483112753\dots$. With a numerical study
we see that the argument of $\Lambda(\tfrac12+it)$ goes from $\pi$ to
a value near $\pi/2$ when $0\le t\le 10$.

\begin{theorem}
There is a  real function $u(t)$ such that $u(t)>0$ for $t>0$ and
such that
\begin{displaymath}
Z(t)=u(t) \Bigl\{-\frac{\pi}{2}+\alpha t+\arctan
2t-\sum_{n\in\Z}\arg\Bigl(1-\frac{s}{b_n}\Bigr)+\Im\frac{s}{b_n}\Bigr\},\qquad
s=\tfrac12+it.
\end{displaymath}
\end{theorem}

\begin{proof}
By definition of $\vartheta(t)$ we have $\pi^{-s/2} \Gamma(s/2)=f(t)
e^{i\vartheta(t)}$, for $s=\tfrac12+it$, where $f(t)=|\pi^{-s/2}
\Gamma(s/2)|>0$. Hence by \eqref{E:ZL}
\begin{align*}
Z(t)=&2\Re\bigl\{e^{i\vartheta(t)}\Lzeta(\tfrac12+it)\bigr\}=
\frac{2}{f(t)}\Re\bigl\{\pi^{-s/2}\Gamma(s/2)\Lzeta(s)\bigr\}\\
&= \frac{2}{f(t)}|\Lambda(\tfrac12+it)|\cos\bigl(\arg
\Lambda(\tfrac12+it)\bigr)
\end{align*}
Define now $a(t)=\frac{\pi}{2}-\arg \Lambda(\tfrac12+it)$, then
\begin{align*}
Z(t)&=\frac{2|\Lambda(\tfrac12+it)|}{f(t)}\sin a(t)\\
&=\frac{2|\Lambda(\tfrac12+it)|}{f(t)}\frac{\sin a(t)}{a(t)}\cdot a(t)
\end{align*}
By Proposition \ref{P:argf} for $t>10$ we have $|a(t)|<0.19$. 
Hence $\sin a(t)/a(t)  >0$.  Thus, we have $Z(t)=u(t) \cdot
a(t)$ were $u(t)>0$  for $t>10$.

Now by Theorem \ref{T:product} we have
\begin{equation}\label{E:arg}
\arg\Lambda(\tfrac12+it)=\pi-\alpha t-\arctan
2t+\sum_{n\in\Z}\Bigl\{\Im\frac{s}{b_n}+\arg\Bigl(1-\frac{s}{b_n}\Bigr)\Bigr\},\qquad
s=\tfrac12+it.
\end{equation}
We know that $1-s/b_n$ is never a negative number for $s=\tfrac12+it$
with $t$ real. In fact if $s$ is such that $1-s/b_n=-x$ we will have
$s=(1+x)b_n$ and these numbers always have a real part greater than
$11$.  Hence, $\arg(1-s/b_n)$ is well defined if we take $\arg z$ such
that $|\arg z|<\pi$. Taking these values of the argument, the series
in \eqref{E:arg} is convergent.

The value of the argument given by \eqref{E:arg} coincides with that
considered in Proposition \ref{P:argf}. To see this, we must only
check a value of $s$, since the two are continuous determination of
the argument. If we take $s=\tfrac12$ we know that the argument
considered in Proposition \ref{P:argf} takes the value $\pi$ on
$s=\tfrac12$.   Now the series for $s=\tfrac12$ that is $t=0$ takes the
value $\pi$ since the term corresponding to $n=0$ is zero and the
one for $n\ne 0$ can be paired with that of $-n$, that
sums to zero.

To obtain our expression for $Z(t)$, we need only to recall the
definition of $a(t)$.

Now for $t>10$ by Proposition \ref{P:argf} the function $a(t)$ has an
absolute value less than $0.19$. Hence the sinus can only vanish
when $a$ vanishes.  Since all zeros of $Z(t)$ for $t>0$ are in
fact $t>10$, we get that $a(\gamma)=0$ for every $\gamma>0$ such
that $Z(\gamma)=0$.  It is obvious that a zero $\gamma$ of
multiplicity $n$ of $Z(t)$ must  also be a zero of the same
multiplicity for $a(t)$. Hence, $Z(t)/a(t)$ is for $t>10$ a
non-vanishing continuous  function.

It is clear that $u(t)$ does not vanish for $0<t<10$ because $Z(t)$
does not vanish, and it is positive in this range since it is
continuous.
\end{proof}

\begin{remark}
The function $u(t)$ and $a(t)$ can be defined for $t<0$ but the
function $u(t)$ is not positive and $a(t)$ does not vanish for $t<0$
at the points where $Z(t)=0$.
\end{remark}

Thus, we consider the function
\begin{displaymath}
a(t)=\frac{\pi}{2}-\arg\Lambda(\tfrac12+it)=-\frac{\pi}{2}+\alpha
t+\arctan
2t-\sum_{n\in\Z}\arg\Bigl(1-\frac{s}{b_n}\Bigr)+\Im\frac{s}{b_n}.
\end{displaymath}
We have, for $t>10$
\begin{displaymath}
|a(t)|=\frac{f(t)}{2|\Lambda(\tfrac12+it)|}\frac{a(t)}{\sin
a(t)}|Z(t)|.
\end{displaymath}
For $t\to+\infty$ we know that
\begin{displaymath}
f(t)\asymp t^{-\frac{1}{4}}e^{-\frac{\pi t}{4}};\quad
|\Lambda(\tfrac12+it)|\asymp t^{-1};\quad \frac{a(t)}{\sin a(t)}\asymp
1.
\end{displaymath}
Thus,
\begin{displaymath}
|a(t)|\asymp  t^{\frac{3}{4}}e^{-\frac{\pi t}{4}}|Z(t)|.
\end{displaymath}

\begin{remark}
It follows that the series for $a(t)$ is useless in locating the zeros
of $Z(t)$ by computation. This series is slowly convergent. The
general term, for a given value of $t$,  is of order $\Im(s/b_n)^2$.
Since $|b_n|\asymp n/\log n$, we need an exponential number of terms
to achieve terms with the order of the sum. In addition, we need to
compute with $\approx\frac{\pi t}{4}$ digits to obtain
significant results. Nevertheless, the expression can have other
theoretical applications.
\end{remark}

\end{document}